\magnification=1200
\hsize=6truein
\hoffset=.5truein
\vsize=8truein 
\parskip=3pt
\tolerance= 10000
\font\rmtwelve=cmbx10 at 12pt
\vskip 3pt

\centerline{\rmtwelve 
Symplectic genus, minimal genus and diffeomorphisms}
\vskip 40pt
\bigbreak
\centerline{Bang-He Li and Tian-Jun Li}
\bigbreak
\bigbreak

\noindent{\bf Abstract}. In this paper,  the symplectic genus for any $2-$dimensional 
class in a $4-$manifold admitting a symplectic structure is introduced, and  its
relation with the minimal genus is studied.
  It is used to describe
which classes in rational and irrational ruled manifolds are realized by connected symplectic surfaces.
In particular, we completely determine which classes with square at least $-1$ in such manifolds can be represented by
embedded spheres. Moreover, based on a new characterization of the action of the diffeomorphism group 
on the intersection forms of a rational manifold, we are able to determine the orbits
of the diffeomorphism group on the set of classes represented by embedded spheres of square at least $-1$
in any  $4-$manifold admitting a symplectic structure.

\bigbreak
\bigbreak

\noindent{\bf \S1 Introduction}
\medbreak

Let $M$ be a smooth,  closed oriented $4-$manifold. 
An orientation-compatible 
symplectic form
on $M$ is a closed two$-$form $\omega$ such that $\omega\wedge \omega$
is nowhere vanishing and agrees with the orientation.
For any
oriented $4-$manifold $M$, its symplectic cone ${\cal C}_M$ is defined as
the set of cohomology classes which are represented by 
orientation-compatible
symplectic forms.

For any class $e\in H^2(M;{\bf Z})$, its minimal genus $m(e)$
is the minimal genus of  a smoothly
embedded connected surface representing the Poincar\'e dual   PD$(e)$.
The problem of determining the minimal genus has involved 
many of the important techniques in $4-$manifold topology, and it
bears its origin in 
the older problem of representing the Poincar\'e dual to a class by 
an embedded sphere (See the 
 excellent survey papers [La1-2] and [Kr1] on these two problems).   

We are here interested in studying both these problems for  $4-$manifolds
with non-empty symplectic cone.     
We will introduce 
the notion of  symplectic genus $\eta(e)$ for $4-$manifolds
with non-empty symplectic cone. 
Recall that any symplectic structure $\omega$
determines a homotopy class of compatible almost complex structures
on the cotangent bundle,
whose first Chern class is called the  canonical class of $\omega$. 
Roughly, the
  symplectic genus $\eta(e)$ of a class $e$ is given by the formula
  $[e^2 + K\cdot e]/2 + 1$, where $K$ has largest pairing against $e$ amongst canonical
  classes of symplectic structures for which the symplectic area of $e$ is
  positive.

$\eta(e)$ has many nice properties,
among which are (i)  invariance under the action of diffeomorphism 
group and (ii) bounding the minimal genus from below. 
We  speculate that, for most  class of positive square, the symplectic
genus is in fact the minimal genus, at least when $b^+(M)=1$ ($b^+(M)$
is the maximal dimension of a positive definite subspace 
of $H^2(M;{\bf R})$).
The minimal genus, by definition, is non-negative. 
And it is easy to see that the  symplectic genus  of a 
sufficiently large multiple of a class of positive square is positive. 
However it is not obvious that the
symplectic genus of any class of positive square is non-negative. 
In this paper we  prove   

\noindent{\bf Theorem A}. Let $M$ be a smooth, closed  oriented
$4-$manifold with non-empty symplectic cone and 
$b^+(M)=1$. Then the symplectic genus of any class of positive square is
non-negative, and it coincides with the minimal
genus for  any sufficiently large multiple of such a class.

The proof of Theorem A is not very difficult except when the manifold
is  a non-minimal rational or irrational ruled manifold. 
In fact, for this kind of manifold
 we are able to obtain a much stronger result. 
Let us explain what such a manifold is. 
Let ${\cal E}_M$ be the set of integral cohomology classes
whose Poincar\'e duals are  represented by
smoothly embedded spheres of squares $-1$. $M$ is said to be (smoothly)
minimal if ${\cal E}_M$ is the empty set.  Any manifold $M$ can be
decomposed as a connected sum of a minimal manifold $N$  with some
number of $\overline {CP}^2$. Such a decomposition is called
a (smooth) minimal reduction of $M$, and $N$ is
a minimal model
of $M$. $M$ is
said to be 
 rational if one of its minimal models is $CP^2$ or $S^2\times S^2$; 
and irrational ruled if one of its minimal models is an
 $S^2-$bundle over a Riemann surface of positive genus. 
When $M$ has non-empty symplectic cone and 
is not rational or irrational ruled, 
$M$ has a unique minimal reduction (see [L1] and [Mc3]).
Using the invariance of the symplectic genus under diffeomorphisms
and the Taubes-Seiberg-Witten theory, we are able to show

\noindent{\bf Theorem B}. Let $M$ be a rational or irrational ruled 
$4-$manifold.  
If $e$ is a class with square at least one, then 
its symplectic genus is non-negative and 
computable. 
Furthermore, if $e\cdot e\geq \eta(e)-1$, then PD$(e)$ is represented
by a connected symplectic surface, and therefore its minimal 
genus coincides with its symplectic genus.

For classes with square zero and $-1$ on rational and irrational
ruled manifolds,
we have similar results. 

Observe that if PD$(e)$ is represented by an embedded sphere, then
$m(e)=0$ and 
therefore $\eta(e)$ is zero as well. 
It turns out that this simple fact enables us to
 completely determine which class of square at least $-1$
is represented by a smoothly embedded sphere in 
any symplectic $4-$manifold.
 When $M$ has nonempty ${\cal C}_M$ and is not rational or irrational ruled,
such a description is known (see [T2], [Mc3] and [L1]). Let 
 $N\# n \overline {CP^2}$ be the unique minimal reduction of
$M$, then, if $e$ has square at least $-1$,    PD$(e)$
is represented by a smoothly  embedded sphere  if and only
if $e$
is a generator of one of the $\overline{CP}^2$.
For rational and ruled manifolds, we have

\noindent{\bf Theorem C}. Let  $M$ be a  rational or irrational ruled manifold and
$e\in H^2(M)$ be a class with square at least $-1$. 
  If
$\eta(e)=0$, then  
PD$(e)$ is represented by a smoothly embedded sphere.
Furthermore, if
PD$(e)$ is represented by a smoothly embedded sphere,
then either $\eta(e)=0$ or $e$ is a non-primitive class of square
zero with $e=pe'$ and $\eta(e')=0$.

We would like to remark that the proofs of Theorems A, B and C are built out
of the work Taubes on Seiberg-Witten invariants realizing symplectic surfaces, 
the wall crossing formula for proving the non-triviality of the Seiberg-Witten invariants,
and the fact that for minimal manifolds with $b^+=1$ it is easy to force symplectic 
surfaces to be connected.
For non-minimal manifolds, we need the additional technical notion of the reduced class.

Beyond determining the set of classes represented by spheres and with square at least
 $-1$, we are also able to
determine the action of the diffeomorphism group on this set. 
 Let us call a class spherically representable if its Poincar\'e dual
is represented
by a smoothly embedded sphere.
Let ${\cal SPH}(M)$ be the set of spherically representable classes
and ${\cal SPH}_{\geq -1}$ be the subset of classes with square 
at least $-1$.
Obviously Diff$(M)$ acts on ${\cal SPH}(M)$ and preserves 
${\cal SPH}_{\geq -1}$.
We are able to completely determine the orbits of ${\cal SPH}_{\geq -1}$ 
 under Diff$(M)$. To state the result we need to introduce more notations. 
We say a class $e$ is of divisibility $p$ if $e=p\tilde e$ with
$\tilde e$ a primitive class.
Let ${\cal SPH}_{s,p}(M)$ be the subset of classes in  
${\cal SPH}(M)$ with
square $s$ and
divisibility $p$. ${\cal SPH}_{s,p}(M)$ can be further 
decomposed depending on the type, i.e. 
whether a class is characteristic or ordinary. 
Recall a class is called characteristic if 
it is an integral lift of the second Stiefel-Whitney class.
Such a class $v$ satisfies 
$v\cdot u=u\cdot u\pmod 2$ for any class $u$. A class is called 
ordinary if it is not characteristic.
Define ${\cal
SPH}^o_{s,p}(M)$ and ${\cal SPH}^c_{s,p}(M)$  
to be the subsets of
 ordinary and characteristic classes in ${\cal SPH}_{s,p}(M)$. 
When 
the group of diffeomorphisms Diff$(M)$ acts on $H^2(M)$, it preserves the
square, the divisibility and the 
type.   Therefore, Diff$(M)$ acts on ${\cal
SPH}^o_{s,p}(M)$ and ${\cal SPH}^c_{s,p}(M)$ separately. 
Remarkably this action is transitive if $s\geq -1$. 

\noindent{\bf Theorem D}.  Let $M$ be a smooth, closed oriented
$4-$manifold with ${\cal C}_M$ nonempty.  Then 
Diff$(M)$ acts transitively on ${\cal SPH}^o_{s,p}(M)$ 
and 
${\cal SPH}^c_{s,p}$ 
when $s\geq -1$.  

The difficult case in Theorem D is when $M$ is a rational manifold
and $s\geq 0$.  
The proof in this case relies crucially on a new  characterization of
the action   of Diff$(M)$ on $H^2(M;{\bf Z})$ in terms of 
the $K-$symplectic cones.

We do not know whether the transitivity continues to hold when
$s$ is less than $-1$.  The case  
$s=-2$ is particularly interesting and will be the subject of further investigation.
We  remark that 
 Theorem D will be applied  in [L2] to prove
that  the fiber sums of relatively minimal Lefschetz fibrations
are minimal manifolds.

\noindent{\bf Convention}. 
From now on, when we say an integral cohomology class
is represented by a surface, we mean  its
Poincar\'e dual is represented by a surface.

The organization of the paper is as follows. In \S2,
we study which automorphisms of the cohomology lattice of a rational
manifold are realized by diffeomorphisms.
Based on a  characterization  by Friedman and Morgan in [FM1-2], 
we give a new  characterization in terms of the $K-$symplectic cones.
This new characterization will be used in \S4 and  is one major new theoretical input in this paper. 
In \S3, we systematically study the
symplectic genus and prove Theorem A and B. Most of this section is
a series of computational lemmas which give enough case-by-case control to
prove the theorems. In \S4, we 
study the problem of representing a class by spheres and 
determine
 the action of 
diffeomorphism groups on ${\cal SPH}(M)$. Theorems C and D will be 
proved there. 

  The authors would like to thank   Janos Koll\'ar,
Ronnie  Lee, Robert Friedman and Gang Tian for helpful discussions. 
This research is partially supported by NSF and 973 Program of P. R. China.

\bigskip
\noindent{\bf \S2. Diffeomorphism group of rational and $K-$symplectic cones}
\medskip

 On a manifold $M$, each
diffeomorphism induces an automorphism
of the lattice of the second integral cohomology.
Hence there is a natural map from 
Diff$(M)$ to the automorphism group of the lattice. Let D$(M)$
be the  image of this natural map. In other words, an automorphism is in D$(M)$
if it is realized by an orientation-preserving diffeomorphism. We will 
describe D$(M)$ for both
rational and irrational ruled manifolds. 

On each irrational ruled manifold $M$, there is a  class
(unique up to sign)  with square
zero 
whose Poinc\'are dual is represented by an embedded sphere. 
It is proved in [FM2] that an automorphism of the cohomology lattice
is in D$(M)$ if and only if that class is preserved up to sign. 
In particular, the -Id automorphism is in D$(M)$.

The case for rational manifolds is rather complicated.
Each rational manifold $M$ is of the form $CP^2\# n {\overline {CP}}^2$.  
When $n\leq 9$,  a result of Wall states that any automorphism is realized by
a diffeomorphism.
The more difficult case $n\geq 10$ requires 
 the concepts of $P-$cell and super 
$P-$cell introduced by Friedman and Morgan [FM1], and  a 
characterization of the 
diffeomorphism group
via these terms. In fact they are not adequate for the purpose
of this paper, and 
we need to consider their partial compactifications and relate them
to the $K-$symplectic cones.

Suppose $M$ is an oriented closed manifold with $b^+=1$, $b^-=n$
and no torsion in $H^2(M;{\bf Z})$. A basis 
$(x, \alpha_1, ..., \alpha_n)$ for $H^2(M;{\bf Z})$ is called 
standard if $x^2=1$, and $\alpha_i^2=-1$ for each $i=1, ..., n$. 
Let $$\eqalign {{\cal P}&=\{e\in H^2(M;{\bf R})|e\cdot e>0\}\cr
{\cal B}&=\{e\in H^2(M;{\bf R})|e\cdot e=0\}\cr
{\overline {\cal P}}&=\{e\in H^2(M;{\bf R})|e\cdot e\geq 0\}.\cr}$$ 
For each class $x\in H^2(M;{\bf Z})$ with $x^2<0$, 
we define $x^{\perp}\in H^2(M;{\bf R})$ to be the orthogonal 
subspace to $x$ with respect to the cup product, and we call
$(x^{\perp})\cap {\cal P}$ the wall in ${\cal P}$ defined by
$x$. 
Let ${\cal W}_1$ be the set of walls in ${\cal P}$ defined by 
integral classes with square $-1$. A chamber for ${\cal W}_1$
is the closure in ${\cal P}$ of a connected component of 
${\cal P}-\cup_{W\in {\cal W}_1}W.$

Any point $x\in {\cal P}$ with square 1  at which  $n$ mutually perpendicular
walls of ${\cal W}_1$ meet is called a corner. 
Any corner is an integral class (see Lemma 2.2
in [FM1]). 
Suppose $C$ is a chamber for ${\cal W}_1$. If $x$ is a corner
in $C$, a standard basis $(x, \alpha_1,\dots, \alpha_n)$ 
for $H^2(M;{\bf Z})$ is called a standard basis adapted to $C$ if
$\alpha_i\cdot C\geq 0$ for each $i$. The canonical class 
of the pair $(x, C)$ is defined
to be $\kappa(x, C)=3x-\sum_i \alpha_i$. 
Suppose $C$ is a chamber for ${\cal W}_1$
 and $x$ is a corner in $C$,  we define
$$P(x, C)=C\cap \{e\in {\cal P}|\kappa(x, C)\cdot e\geq 0\}.$$
Any subset of ${\cal P}$ of the form $P(x, C)$ is called a $P-$cell. 

\noindent {\bf Notation}. For any $U\subset {\cal P}$ (similarly ${\cal B}$, $\overline{\cal P}$), 
 we will use int$_{\cal P}(U)$ (similarly int$_{\cal B}(U)$, int$_{\overline {\cal P}}(U)$) to denote 
$U\cap$int$({\cal P})$ (similarly $U\cap$int$({\cal B})$, 
$U\cap$int$(\overline {\cal P})$). For any $V\subset \overline{\cal P}$ (respectively ${\cal P}$),
we will use $\bar V$ (similarly $\hat V$) to denote its closure in $\overline{\cal P}$ (similarly
${\cal P}$).

The basic properties of $P-$cells are summarized in the following lemma.

\noindent{\bf Lemma 2.1}. 

\noindent 1.
A $P-$cell is a chamber for the set of walls 
${\cal W}_1 \cup \{\kappa(x,C)^{\perp}\cap {\cal P}\}$.

\noindent 2. If $P(x, C)=P(x', C')$, then $\kappa(x, C)=\kappa(x', C')$.
Thus for each $P-$cell $P$ we can assign  a unique canonical
class of the form
$\kappa(x, C)$, which  will be written as
$\kappa(P)$.  

\noindent 3. If $\psi$ is an automorphism of the lattice and 
$P$ is a $P-$cell with canonical class $\kappa$, then
$\psi\cdot P$ is also a $P-$cell with canonical class $\psi(\kappa)$.

\noindent 4. If $P$ and $P'$ are distinct $P-$cells, then
int$_{\cal P}(P)$ and int$_{\cal P}(P')$ are disjoint. 

\noindent 5. If $P$ and $P'$ are distinct $P-$cells, then
int$_{\cal B}(\overline P\cap {\cal B})$ and int$_{\cal B}(\overline P'\cap {\cal B})$ 
are disjoint. In other words,
the interiors of the ${\cal B}-$boundaries of the closure of 
distinct  $P-$cells are also disjoint.

\noindent{\it Proof}. The proofs of the  first $4$ properties 
can be found  in chapter II in [FM 1].
Here we prove property $5$. 
Notice that $\overline {\cal P}={\cal P}\cup {\cal B}$.  
If $x$ is  any point in int$_{\cal B}(\overline P\cap {\cal B})$, then
the intersection of ${\cal P}$ with any sufficiently small neighborhood
 in $\overline {\cal P}$ 
 of $x$ is non-empty and is contained in int$_{\cal P}(P)$. 
Thus if
int$_{\cal B}(\overline P\cap {\cal B})$ and int$_{\cal B}(\overline P'\cap {\cal B})$ 
intersect, 
then int$_{\cal P}(P)$ and int$_{\cal P}(P')$ overlap and hence they are the same
$P-$cells by
property $4$. Therefore distinct $P-$cells have disjoint ${\cal B}-$boundaries.

It turns out that $P-$cells are closely related to the 
$K-$symplectic cone introduced in [LLiu1]. Let us  recall the definition
of $K-$symplectic cone. A class $K\in H^2(M;{\bf Z})$ is called a symplectic canonical
class if it is the canonical class of some orientation-compatible symplectic structures.
Let ${\cal K}$ be the set of symplectic canonical classes.
For any $K\in{\cal K}$ we introduce the $K-$symplectic cone:  
$${\cal C}_{K}=\{e\in {\cal C}_M| \hbox{ $e=[\omega]$ for some $\omega \in \Omega_K $ }\},$$
where $\Omega_K$ is the set of orientation-compatible 
symplectic forms with $K$ as the symplectic canonical class. 
It is shown in [LLiu2] that  
${\cal C}_{K}$
is disjoint from ${\cal C}_{K'}$ if $K\ne K'$. 
For a manifold with $b^+=1$ and any $K\in {\cal K}$, we can in fact determine
${\cal C}_K$ in terms of a certain subset of ${\cal E}_M$. Recall that
 ${\cal E}_M$ is the set of integral cohomology classes
represented by
smoothly embedded spheres of square $-1$. When there is no confusion
we will omit the subscript $M$.   
 Introduce the set of $K-$exceptional spheres as
$${\cal E}_K=\{E\in {\cal E}|K\cdot E=-1\}.$$
It is proved in Theorem 4 in [LLiu1] that 
$${\cal C}_K=\{e\in {\cal P}|e\cdot E>0 \hbox{ for any $E\in {\cal E}_K$ }\}.$$
Let ${\hat {\cal C}_K}$ be the closure of
${\cal C}_K$ in ${\cal P}$. Then it is not hard to prove 
 $${\hat {\cal C}_K}=\{e\in {\cal P}|e\cdot E\geq 0 
\hbox{ for any $E\in {\cal E}_K$ }\}.$$

 In order to link the $P-$cells and the symplectic cones,  
 we also need to consider
 good generic surfaces
as in [FM1].
A good generic  surface $X$ is an algebraic  surface such that 
the anti-canonical divisor is effective and smooth, and that 
any smooth rational curve has square no less than $-1$.
All such  surfaces are rational surfaces
and can be holomorphically blown down to $CP^2$ (see I.2 in [FM1]). 
Let $\rho:X\longrightarrow CP^2$ be a holomorphic blow down
with exceptional fibers $F_1, ..., F_n$, where each
$F_i$ is an embedded rational curve with square $-1$. 
Let $K_X$ be the canonical class
of $X$. Then $-K_X=3\rho^*(H)-\sum_{i=1}^n F_i$, where
$H$ is a hyperplane section of $CP^2$. 

The surface $X$ has many K\"ahler metrics. Associated with each such metric
is its K\"ahler form and associated cohomology class in $H^2(X;{\bf R})$.
The K\"ahler cone ${\cal A}(X)$ of $X$  is then 
 the set of all K\"ahler cohomology classes.  
 By the Nakai-Moishezon criterion, 
the K\"ahler cone ${\cal A}(X)$ consists of all the classes 
in ${\cal P}$ which pair positively on
 any holomorphic curve. 
Let $\hat {\cal A}(X)$ be the closure of ${\cal A}(X)$ in
${\cal P}$.

 
\noindent{\bf Proposition 2.2}. Let $X$ be a good generic surface. 
Let $P_0$ be the $P-$cell containing the class $\rho^*H$. Then 
$P_0$ coincides with $\hat{\cal A}(X)$, 
  and $\kappa(P_0)=-K_X$. Moreover, 
$$P_0=\{e\in {\hat {\cal C}_{K_{X}}}|e\cdot (- K_{X})\geq 0\}.$$

\noindent{\it Proof}. The first statement is proved in 
II. 3 and II.4 in [FM1]. So we only need to show that $\hat{\cal A}(X)$ 
consists of all classes in  $\hat{\cal C}_{K_X}$ which pair non-negatively 
with $(-K_X)$.  

Since a K\"ahler form is a symplectic form, the K\"ahler cone
${\cal A}(X)$
is certainly a subset of
the $K_X-$symplectic cone ${\cal C}_{K_X}$.  Therefore 
$\hat {\cal A}(X)$ is a subset of ${\hat {\cal C}_{K_{X}}}$.
To prove the inclusion 
in the  other direction, we need the following
result:
$$\hat {\cal A}(X)=\{e\in {\cal P}|e\cdot (-K)\geq 0
\hbox { and }  e\cdot E\geq 0 
\hbox{ for any $E\in {\cal E}^{hol}(X)$ }\}, $$
which is  Proposition 3.4 in [FM1].
Here ${\cal E}^{hol}(X)$ is the set of embedded rational curves with
square $-1$. 
 With this characterization of $\hat {\cal A}(X)$ we just have to show that,
 on a good generic  surface, 
any class $e\in {\hat {\cal C}_{K_X}}$ with $e\cdot (-K_X)\geq 0$ is 
non-negative on any class in ${\cal E}^{hol}(X)$.
This follows from the obvious inclusion ${\cal E}^{hol}(X)\subset 
{\cal E}_{K_X}$.
The proof is complete.

\noindent{\bf Remark 2.3}. From Proposition 2.2, we find
$$P_0=\{e\in {\cal P}|e\cdot (- K_{X})\geq 0 \hbox { and }  e\cdot 
E\geq 0
\hbox{ for any $E\in {\cal E}^{hol}(X)$ }\}.$$ 
Since $P_0$ coincides with ${\hat A}(X)$, it is possible that the two sets
${\cal E}^{hol}(X)$
and ${\cal E}_{K_X}$ are the same.                         

\noindent {\bf Lemma 2.4}. Let $M$ be a rational $4-$manifold. 
For each $K\in {\cal K}$, there exists a $P-$cell $P_K$
such that 
$\kappa(P_K)=-K$ and 
$$P_K=\{e\in {\hat {\cal C}_K}|e\cdot (-K)\geq 0\}.$$

\noindent{\it Proof}. Suppose $X$ is a good generic surface
and $M$ is the underlying rational $4-$manifold. 
By the result in [LLiu1] that
 Diff$(M)$ acts transitively on ${\cal K}$, there is a diffeomorphism
$\phi$ of $M$ such that $\phi^*(K_X)=K$. 
Since $\phi^*\omega\in { \Omega}_{\phi^*(K_X)}$ for any 
$\omega\in \Omega_{K_X}$, we have
  ${\cal C}_K=\phi^*{\cal C}_{K_X}$. Thus by Proposition 2.2 we have 
$$\phi^*(P_0)=\{e\in {\hat {\cal C}_{K}}|e\cdot (- K)\geq 0\}.$$
Let $P_K=\phi^* P_0$. By Lemma 2.1(2-3), $P_K$ is still a $P-$cell
with canonical class $K$. 
 The proof is complete. 


Now we introduce  super $P-$cell, which is defined via
a reflection group associated to a $P-$cell. 
Suppose $\gamma$ is a class with square $-1$ or $-2$. We can define
an automorphism of the lattice as follows, 
$$R(\gamma)\beta=\beta+  {2(\gamma\cdot \beta)\over {|\gamma\cdot \gamma|}}
\gamma.$$
This automorphism $R(\gamma)$ is called the reflection along $\gamma$.
For a $P-$cell $P$ define ${\cal G}_P$ to be the set 
$$\{\alpha|\alpha^2=-1, \alpha\ne \kappa(P) \hbox{ and  $\alpha$ defines
a wall of ${\cal P}$}\}.$$
Let ${\cal R}(P)$ be the group generated by reflections along classes in
${\cal G}_P$.  

The super $P-$cell of $P$  is defined as
$$S(P)=\cup_{\psi\in {\cal R}(P)} \psi( P).$$

We will need the following simple fact on reflections.

\noindent{\bf Lemma 2.5}. Suppose $F=\psi (F_0)$ for some
$F_0\in {\cal G}_P$. Then $R(F)=\psi\circ R(F_0)\circ \psi^{-1}$.
In particular, $R(F)\in {\cal R}(P)$.

\noindent {\it Proof}. For any class $x$, we have
$$\eqalign {R(F)\circ \psi (x)
&=\psi(x)+ 2(F\cdot \psi(x))F=\psi(x)+2(\psi(F_0)\cdot \psi(x))\psi(F_0)\cr
\psi\circ R(F_0)(x)&=\psi(x+2(F_0\cdot x)F_0)=\psi(x)+2(\psi(F_0)\cdot \psi(x)\psi(F_0).\cr}$$ 
 So $R(F)\circ \psi=\psi\circ R(F_0),$ and the 
statements follow. 

\noindent {\bf Proposition 2.6}. Let $M$ be a rational $4-$manifold.
Any good generic surface $X$ with $M$ as the underlying manifold
gives rise to a $P-$cell of $M$, denoted by $P_0$.

\noindent 1. If $\phi$ is an automorphism, then 
$\phi( S(P))$ is also a super $P-$cell. In particular,
$-S(P)$ is a super $P-$cell. 

\noindent 2. An automorphism is in D$(M)$
 if and only if it preserves the distinguished 
super $P-$cell $S(P_0)$ up to sign. Consequently, D$(M)$ is generated by -Id, ${\cal R}(P_0)$ and
the isotropy subgroup of $P_0$.

\noindent 3. If int$_{\cal P}(S(P))\cap$ int$_{\cal P}( S(P'))$ is not empty, then $S(P)=S(P')$.

\noindent 4. If int$_{\cal B}(\overline {S(P)}\cap {\cal B})$ and 
int$_{\cal B}(\overline {S(P')}\cap {\cal B})$ intersect, then $S(P)=S(P')$.

\noindent 5. If int$_{\overline {\cal P}}(   \overline {\phi( S(P))} )$ 
and int$_{\overline {\cal P}}(   \overline {\phi( S(P'))} )$
intersect, then $S(P)=S(P')$.

\noindent{\it Proof}. The first three properties are in chapter II in [FM1]. 
 The proof of property 4 goes exactly along the line of the proof of
the analogous property for the $P-$cells in Lemma 2.1. 
If $x$ is  any point in int$_{\cal B}(\overline {S(P)}\cap {\cal B})$, then
the intersection of ${\cal P}$ with any sufficiently small neighborhood
 in $\overline {\cal P}$ 
 of $x$ is non-empty and is contained in int$_{\cal P}(S(P))$.  
Thus if  int$_{\cal B}(\overline {S(P)}\cap {\cal B})$ and 
int$_{\cal B}(\overline {S(P')}\cap {\cal B})$ 
intersect, 
then $S(P)$ and $S(P')$  have overlapping interiors and hence they are the same
super $P-$cells by
property 3.

Notice that 
$$\hbox{int}_{\overline {\cal P}}(  \overline {\phi( S(P))} )=
\hbox{int}_{\cal P}
(S(P))\cup
\hbox{int}_{\cal B}(\overline {S(P)}\cap {\cal B}).$$
The last statement follows immediately from the properties 3 and 4. 
The proof is complete. 

In the next proposition we are going to relate the
super $P-$cells $\pm S(P_0)$ to the $K-$symplectic cones. 

\noindent{\bf Proposition 2.7}. Define $M$ and $X$  as in Proposition 2.5. 
 Let 
$K_0$ be the canonical class of $P_0$. Then
every $K\in{\cal K}$ 
is of the form $\pm \psi (K_0)$, where $\psi\in {\cal R}(P_0)$.
Consequently, 
$$\eqalign {S(P_0)\cup -S(P_0)&=\cup_{K\in {\cal K}} P_K\cr
\overline{S(P_0)\cup -S(P_0)}\cap {\cal B}&=\cup_{K\in {\cal K}}
\overline{ P}_K\cap {\cal B}.\cr}$$

\noindent{\it Proof}. 
This is a consequence of the result in [LLiu1] which states
that D$(M)$ acts transitively on ${\cal K}$.  The positive cone ${\cal P}$ has two connected 
components.
Let $K$ be a symplectic canonical class such that
${\cal C}_K$ and ${\cal C}_{K_0}$ are in the same connected component of
${\cal P}$. Since D$(M)$ acts transitively on ${\cal K}$, there
exists $\psi'\in$D$(M)$ such that $\psi(K_0)=K$. By
Proposition 2.6(2), $\psi'(S(P_0))=\pm S(P_0)$.  Since $\psi'(P_0)$
is still in the same component of $P_0$, $\psi'(S(P_0))= 
S(P_0)$. Therefore $\psi'(P_0)$ is a $P-$cell within $S(P_0)$.
 By the definition of a super $P-$cell,
$\psi'(P_0)=\psi(P_0)$ for some $\psi \in {\cal R}(P)$.
Therefore $K=\psi (K_0)$.  
By Lemma 2.1, $\psi (P_{K_0})$ and $P_K$ have the same
canonical class and therefore they are identical. 
Notice we have shown that 
$$\cup_{K\in {\cal K}} P_K \subset S(P_0)\cup -S(P_0).$$
To prove the inclusion in the other direction, we just need to 
show that, for any $\psi\in {\cal R}(P)$, $\psi(P_0)=P_K$
for some $K\in {\cal K}.$ This is obvious since 
$K=\psi(K_0)$ is certainly a symplectic canonical class.
The proof is finished.  

 It is mentioned in [FM2] that super $P-$cells are chambers for
the walls given by primitive characteristic 
classes with square $9-n$. 
We can in fact show that the set of walls for $S(P_0)$ and $-S(P_0)$ is 
just the 
set of  the symplectic canonical classes.

Now we are able to present the main result of this section,
a characterization of 
D$(M)$ in terms of $K-$symplectic cones.

\noindent {\bf Theorem  2.8}. Let $M$ be a rational $4-$manifold.
 An automorphism $\phi$ is  in D$(M)$ if and only if
there are $K$ and $K'$ in ${\cal K}$ and 

\noindent 1. either there are classes
$e\in \hat {\cal C}_K$ and $e'\in \hat{\cal C}_{K'}$
with
$e\cdot (-K)>0$ and $e'\cdot (-K')>0$, such that
$e$ is mapped to $e'$ by $\phi$,

\noindent 2. or there are classes $e\in \overline {{\cal C}}_K\cap {\cal B} $  
and $e'\in \overline {{\cal C}}_{K'}\cap {\cal B} $ with
$e\cdot (-K)> 0$ and $e'\cdot (-K')> 0$, such that
$e$ is mapped to $e'$ by $\phi$.

\noindent{\it Proof}.  Due to Proposition 2.7, 
in the first case, we just have to show that $e$ and $e'$
are in the interior of $S(P_0)\cup -S(P_0)$. The arguments for
$e$ and $e'$ are exactly the same, so we will only argue
for $e$.  
By Lemma 2.4, $e\in P_K$. 
If $e$ is in the interior of $P_K$, then 
$e$ is in the interior of $S(P_0)\cup -S(P_0)$ by Proposition 2.7.
$e$ may fail to be in the interior of $P_K$ only when
$e\cdot E=0$ for some $E\in {\cal E}_K$. Suppose 
$P_K=\pm \psi(P_0)$ for some $\psi\in {\cal R}(P_0)$, then
$E=\psi (E_0)$ for some $E_0\in {\cal G}_{P_0}$.  By Lemma 2.5,
the $P-$cell obtained by
reflecting $P_K$ along $E$ is still in $S(P_0)\cup -S(P_0)$. Thus we see
that $e$ must be in the interior of $S(P_0)\cup -S(P_0)$. 

The proof in the second case is similar. We just have to show that
$e$ and $e'$ are in $ \hbox{int}_{\cal B}(\overline {S(P_0)\cup
-S(P_0)}\cap {\cal B})$ and we only have to argue for $e$. 
By Lemma 2.4, $e\in  \overline{P}_K\cap {\cal B}$. 
 $e$ may fail to be in the interior of 
$\overline{P}_K\cap {\cal B}$ only when
$e\cdot E=0$ for some $E\in {\cal E}_K$.  However by Lemma 2.5, 
the reflection of $\overline{P}_K \cap {\cal B}$ along $E$ is still in 
$\overline{S(P_0)}\cap {\cal B}$. Thus we see
that $e$ must be in the interior of $\overline{S(P_0)\cup -S(P_0)}\cap {\cal B}$. The 
theorem is proved.

\bigskip
\noindent{\bf \S3. Symplectic genus}

\medskip

We first give the formal definition of the symplectic genus for manifolds with 
non-empty symplectic cone.
For any integral
class $e\in {\cal P}$, 
we first define a subset of ${\cal K}$, 
$${\cal K}_e=\{K\in{\cal K}| \hbox {there exists a class $\tau \in {\cal C}_K$
such that $\tau\cdot e>0$} \}.$$
We further define a subset of ${\cal K}_e$,

$${\cal K}(e)=\{K\in{\cal K}_e| K\cdot e \geq K'\cdot e \hbox { for any $K'\in{\cal K}_e $} \}.$$  

\noindent{\bf Definition 3.1}. Let  
$K$ be any class in $ {\cal K}(e)$. The symplectic genus of $e$ is defined to be
$$\eta(e)={1\over 2}(e\cdot K+e^2)+1.$$

We now list some simple  properties of symplectic genus.

\noindent {\bf Lemma 3.2}.

\noindent 1.  The symplectic genus is no bigger than the minimal genus.
Furthermore, if a class is represented by a connected symplectic surface,
then  its symplectic genus is equal to its minimal genus.  

\noindent 2. $\eta(-e)=\eta(e)$.

\noindent 3. For any positive integer $p$,
$$\eta(pe)=p\eta(e)-(p-1)+{(p^2-p)\over 2}e\cdot e.$$
In particular, $\eta(pe)\ne 0$ when $e\cdot e=0$ and $p\geq 2$. 

\noindent 4. The symplectic genus is 
invariant under the action of the group of orientation-preserving
diffeomorphisms.

\noindent 5. The symplectic genus of any class of a 
sufficiently large multiple of any class of positive square is positive. 

\noindent{\it Proof}.
Property 1 is a consequence of the adjunction inequalities.
When $b^+>1$ the adjunction inequality in  [KM], [MST], [OS] and [T2] 
asserts that
the genus $g$  of any embedded surface representing $e$ satisfies 
$$2g-2\geq |K\cdot e|+e\cdot e  \eqno (3.1)$$
for any symplectic canonical class $K$.

When $b^+=1$ and $e$ has non-negative square, 
the adjunction inequality in [LLiu2] asserts that 
$$2g-2\geq K\cdot e+e\cdot e \eqno (3.2)$$
for any symplectic canonical class $K\in {\cal K}_e$. 

When $e$ has negative square, inequality (3.2) still holds and is  
basically proved  in \S3 in [OS]. We explain here briefly. 
Suppose $\omega$ is a symplectic form whose class $\tau$ pairs positively
with $e$, and let $K(\omega)$ be its symplectical canonical class.
Let $s_0$ be   the canonical Spin$^c$ structure with $c_1(s_0)=-K(\omega)$.
The class $e$ determines another Spin$^c$ structure, denoted by $s_0-e$.
Suppose $e$ is represented by an embedded surface of genus
$h$ such that 
$$2h-2<K(\omega)\cdot e+e\cdot e.$$ Then,  by Theorem 1.3 in [OS] and 
the corresponding result in [FS], 
in a common chamber for both $s_0$ and $s_0-e$ which is perpendicular
to $e$, the Seiberg-Witten invariant of
$s_0$ being nontrivial implies that the invariant 
of $s_0-e$ is nontrivial as well.
The $\omega-$symplectic chamber is such a chamber. Moreover, according
to Taubes ([T1]),
in this chamber, the Seiberg-Witten invariant of $s_0$ is nontrivial.
Therefore, in the $\omega-$symplectic chamber,
 the Seiberg-Witten invariant of $s_0-e$ is nontrivial as well.
By another result of Taubes ([T2]), $\tau\cdot (-e)>0$. 
This  contradicts our assumption, so inequality (3.2) 
still holds in this case.   
Therefore, in any case, we have $m(e)\geq \eta(e)$. 

Suppose that $e$ is represented by a genus 
$h$ symplectic surface with respect to 
a symplectic form $\omega$. Then $\omega$ is positive on this surface. 
If $K(\omega)$ is the symplectic canonical class
of $\omega$, then $K(\omega)\in {\cal K}_e$ and 
$2h-2=K(\omega)\cdot e+e\cdot e.$ Together with inequalities (3.1) and (3.2),
we see that $h=m(e)=\eta(e)$.

If $K$ is the symplectic canonical class of a symplectic form $\omega$, then
$-K$ is the symplectic canonical class of the symplectic form $-\omega$.
Therefore, $${\cal K}_{-e}=\{-K|K\in {\cal K}_e\} \hbox { and }
{\cal K}(-e)=\{-K|K\in {\cal K}(e)\} \eqno (3.3).$$
And $\eta(-e)=\eta(e)$ is an immediate consequence of equation (3.3).

For any positive integer $p$, we have 
$${\cal K}_e={\cal K}_{pe} \hbox { and }  {\cal K}(e)={\cal K}(pe) 
\eqno (3.4).$$
The formula for $\eta(e)$ then follows from equation (3.4) with a straightforward calculation.
When $e\cdot e=0$, $\eta(pe)$ is therefore given by $p(\eta(e)-1)+1$. Evidently it is not divisible by $p$ 
and concequently cannot be zero if
$p\geq 2$. 

 It is shown in Proposition 4.1 in [LLiu1] that,
if $\phi$ is an orientation-preserving diffeomorphism, then
$\phi^*{\cal C}_K={\cal C}_{\phi^*K}.$ 
It follows  that
$$\phi^*{\cal K}_e={\cal K}_{\phi^*e} \quad \hbox{and} \quad \phi^*{\cal K}(e)={\cal K}({\phi^*e}) \eqno (3.5).$$
Property 4, then, is an immediate consequence of equation (3.5).

The last property follows directly from the definition. Let $e$ be a class with positive square. When $N$ is large,
$Ne\cdot Ne$ dominates $Ne\cdot K$ for any $K\in {\cal K}(e)$, and  
therefore $Ne$ has positive symplectic genus. Lemma 3.2 is proved.

Now we set out to prove Theorem B.
The proof requires  the notion of reduced classes
for non-minimal rational and irrational ruled manifolds
(for rational manifolds, it is introduced in [Ki] and [G]).
A nice property of this notion is that every class with positive 
square can be transformed in an explicit way  to  a reduced class via
diffeomorphisms. Thus
 by Lemma 3.2(4) we only have to show that
Theorem B holds for any reduced class $e$. 

To introduce the reduced class let us review the minimal reductions of 
a rational or irrational ruled manifold.
The only minimal rational manifolds are $CP^2$ and $S^2\times S^2$.
 And  a non-minimal rational manifold  has two kinds of decomposition- 
it is either decomposed as $CP^2\# n\overline {CP}^2$ or as 
$S^2\times S^2 \# (n-1)\overline {CP}^2$. 
We will always use the first decomposition and call it a standard 
decomposition.
 The picture for irrational ruled manifolds is
similar. 
$S^2-$bundles over a Riemann surface of positive genus are the only 
minimal irrational ruled manifolds. Fix the base  surface $\Sigma_g$,
there are two $S^2-$bundles over it, the trivial one
$S^2\times \Sigma_g$ and the unique
non-trivial one $S^2\tilde{\times}\Sigma_g$. 
A non-minimal  irrational ruled manifold also has two types of
decomposition,   either  as 
$S^2\times \Sigma_g\# n\overline {CP}^2$ or as 
$ S^2\tilde{\times}\Sigma_g \# n\overline {CP}^2$. We will use the first 
decomposition and call it a standard decomposition.

Let $H$ be a generator of $H^2(CP^2;{\bf Z})$ and $E_1,\dots, E_n$ be 
the generators of the
$\overline{CP}^2$.  
Let $U$ and $T$ be  classes in $S^2\times \Sigma_g$ represented by
$\{pt\}\times \Sigma_g$  and  $S^2\times \{pt\}$ respectively. 
$H, E_1,\dots, E_n$ are naturally considered as classes in
$H^2(CP^2\# n\overline {CP}^2;{\bf Z}) $ and form a basis.
We will call such a basis  a standard basis. Similarly, 
$U, T, E_1, \dots, E_n$ are naturally considered as classes
in  $H^2(S^2\times \Sigma_g\# n\overline {CP}^2; {\bf Z})$ and form a basis.
Such a basis is also called a standard basis. Given such a basis, according to Wall ([W]), 
an automorphism is called trivial if either it permutes the $E_i$ or it is a
  reflection along an $E_i$.  It was  shown in [W] that
trivial automorphisms are in D$(M)$.  

 On $CP^2\# n\overline {CP}^2$, let  $K_{0}=-3H+\sum_i E_i$;
and on $S^2\times \Sigma_g\# n\overline {CP}^2$, let  $K_0=-2U+(2g-2)T+\sum_i E_i$. By the blow up construction (see e.g. [Mc1])
$K_0$ is a symplectic canonical class.

\noindent{\bf Definition 3.3}.
For a non-minimal rational manifold with a standard 
decomposition $CP^2\# n\overline{CP}^2$
and a standard basis $\{H, E_1,\dots, E_n\}$, 
a  class  $\xi=aH-\sum^n_{i=1}b_iE_i$ is called reduced
if $$\cases{b_1\geq b_2\geq \cdots \geq b_n\geq 0\cr
a\geq                b_1+b_2+b_3.\cr}$$
For a non-minimal irrational ruled manifold with a standard decomposition 
$S^2\times \Sigma_g\# n\overline{CP}^2$ and a
standard basis $\{U, T,E_1, \dots, E_n\}$, a class $e=aU+bT-\sum c_iE_i$ is called 
reduced if
$$\cases{a\geq 0, c_1\geq c_2\geq \cdots\geq c_n\geq 0\cr
a\geq c_i
\hbox{ for any }i.\cr}$$

Reduced classes have the following properties:

\noindent{\bf Lemma 3.4}. Let $M$ be a non-minimal rational or 
irrational ruled manifold
with a standard decomposition and a standard basis. 

\noindent 1. Any class of non-negative square is equivalent  to a reduced class
under the action of orientation-preserving diffeomorphisms. 
Moreover we can find such a diffeomorphism by a simple algorithm.

\noindent 2. For a class with square $-1$,  when $b^-(M)\ne  2$, it 
either  has 
reduced form
or   is equivalent to the class $E_1$;  
when $b^-(M)=2$, another possibility is that it is characteristic, and
equivalent to 
$H-E_1-E_2$ in the  rational case and to $T-E_1$ in the 
ruled case. 
Similarly, for a class with square $-2$,  when $b^-(M)\ne  3$,
it  either  has 
reduced form
or   is equivalent to the class $E_1+E_2$;  
when $b^-(M)=3$, another possibility is that it is characteristic, and 
equivalent to 
$H-E_1-E_2-E_3$ in the rational case and to $T-E_1-E_2$ in the 
irrational ruled case. 

\noindent 3.  If $e$ is reduced, then 
 $e\cdot F\geq 0$
for any $F\in {\cal E}_{K_0}$.

\noindent 4. If $e$ is reduced, then $K_0\in {\cal K}_e$.

\noindent 5. 
If $e$ is a reduced class with non-negative square, then $K_0\in{\cal K}(e)$,
and consequently $\eta(e)$ is given by 
$(K_0\cdot e+e\cdot e)/2 +1$.

\noindent{\it Proof}.
We divide the proof into two cases.

\noindent (i). First consider a non-minimal rational manifold  
with a standard decomposition
 $CP^2\# n\overline{CP}^2$
and a standard basis. When $n\leq 9$, all the properties have been 
established (for 1 and 4 see [Li1], for 2 see [LiL2]).
So we assume that $n\geq 10$.  

\noindent {\it Property 1}. In fact, it was
also proved  in [Li1].  Since we will use the similar arguments 
to prove property 2, we provide  some details here.
 Suppose $e=aH-\sum_i b_i E_i$ is a class
with non-negative square.    First of all, by the trivial 
automorphisms, we can arrange so that $a\geq 0$ and
$b_1\geq b_2\geq ...\geq b_n\geq 0$. 
When $n\geq 3$, the class $H-E_1-E_2-E_3$ is represented
by an embedded sphere with square $-2$.  
So the reflection along $H-E_1-E_2-E_3$ is an automorphism in 
D$(M)$. Under this reflection,
 $a$ is mapped to 
$a'=2a-(b_1+b_2+b_3)$. If $e$ is not already reduced 
and $2a-(b_1+b_2+b_3)\geq 0$, then     $0\leq a'<a$. From this we see the
process can be repeated either to lead to a reduced class or to a class with
$2\tilde a-(\tilde b_1+\tilde b_2+\tilde b_3)<0$. 
However, if $2\tilde a<(\tilde b_1+\tilde b_2+\tilde b_3)$, then  from 
$$\tilde a^2\geq \sum_i \tilde b_i^2 \quad \hbox{and} \quad (\tilde
b_1+\tilde b_2+\tilde b_3)^2\leq
3(\tilde b_1^2+\tilde b_2^2+\tilde b_3^2),$$ 
we have
$$(\tilde b_1^2+\tilde b_2^2+\tilde b_3^2) \leq \tilde a^2<
(3/4)(\tilde b_1^2+\tilde b_2^2+\tilde b_3^2).$$
This is an obvious contradiction.

\noindent {\it Property 2.}  
Suppose we have a class $e$ of square $-1$.
The same argument as above proves that $e$ is either equivalent
to a reduced class or a class with $$(\tilde b_1^2+\tilde b_2^2+
\tilde b_3^2)-1 \leq \tilde a^2<
(3/4)(\tilde b_1^2+\tilde b_2^2+\tilde b_3^2).$$
In the latter case $(\tilde b_1^2+\tilde b_2^2+\tilde b_3^2)< 4$ and
 we easily deduce that, up to trivial automorphisms,  
the only such class is    $E_1$.     
For a class with square $-2$,  the same argument again
proves that $e$ is either equivalent
to a reduced class or a class with 
$$(\tilde b_1^2+\tilde b_2^2+\tilde
b_3^2)-2\leq \tilde a^2<
(3/4)(\tilde b_1^2+\tilde b_2^2+\tilde b_3^2).$$
In the latter case we easily find that there are only two possibilities
up to trivial automorphisms: 
$H-E_1-E_2-E_3$   and $E_1+E_2$.    However 
$H-E_1-E_2-E_3$ is equivalent to $E_1+E_2$  under the reflection
along $H-E_2-E_3-E_4$ and trivial automorphisms. 

\noindent{\it Property 3.} 
Assume that 
$F=tH-\sum s_i E_i$. Then 
$$K_0\cdot F=-3t+\sum s_i=-1 \hbox{ and } e\cdot F=at-\sum s_ib_i.$$
 It was shown in
[LLiu2] that $F, E_1,..., E_n$ and $H$ are all represented
by connected smooth $J-$holomorphic spheres for
some almost complex structure $J$. By the positivity of intersection
of distinct $J-$holomorphic curves, $t\geq 0$, and $s_i\geq 0$
unless  $F=E_i$. If $F=E_i$ for some $i$,
clearly we have $e\cdot F\geq 0$. If $F\ne E_i$, then $t> 0$ and
$t\geq s_i\geq 0$.
We can divide the $b_i$ into $t$ groups, each consisting of no more than three
$b_i$. Since $s_i$ is no bigger than 
$t$, the division can be made such that the $b_i$ in each group 
have distinct indices. The condition of $e$ being reduced implies that
$a-b_i-b_j-b_k\geq 0$ for any $i, j, k$ which are mutually distinct.
The property follows.

\noindent{\it Property 4.} 
Notice that for any sufficiently small $\epsilon$, $\omega_{\epsilon}=H-\sum \epsilon E_i$
is a symplectic form with canonical class $K_0$. 
Since $\omega_{\epsilon}\cdot e>0$
for $\epsilon$ small, we have $K_0\in {\cal K}_e$. 

\noindent{\it Property 5.}
Since a reduced class $e$ with non-negative square has a positive $H$ term, by the light cone lemma
in [LLiu2], the class of 
 a symplectic form is positive on 
$e$ only when it has a positive $H$ term as well. 
Therefore, 
if $K$ is  any symplectic canonical class in ${\cal K}_e$, then it has a
negative $H$ term. 
By Proposition 2.5,
any $K\in {\cal K}_e$ is of the form $\psi(K_0)$ for some 
$\psi\in {\cal R}(P_0)$. 
We claim that $K\cdot e\leq K_0\cdot e$ for any
$K$ of the form $K\in {\cal K}_e$. 
Once this is established 
it is clear that $K_0\in {\cal K}(e)$ and property 4 follows. 
Now we proceed to prove the above claim.
Write $\psi$ as
$R^{(F_k)}\circ\dots \circ R^{(F_1)} (K_0)$ 
where $F_i\in {\cal E}_{K_0}$. 
 Since $K_0\cdot F_i=-1$, 
 we have 
$$\eqalign{&R^{(F_1)}(K_0)=K_0+2(K_0\cdot F_1)F_1=K_0-2F_1\cr
&R^{(F_2)}(K_0-F_1)=[K_0-F_1]+2[(K_0-F_1)\cdot F_2]F_2=
K_0-F_1-(1+F_1\cdot F_2)F_2.\cr}$$
This, together with the fact $F_i\cdot F_j\geq 0$
due to positivity of intersection, we have
$$K=K_0-2F_1-\sum_{i=2}^k 2c_iF_i, \quad c_i\geq 1.$$
Now the claim follows from property 3.

\noindent (ii). Now consider a non-minimal
irrational ruled manifold with a standard decomposition 
and a standard basis.
Suppose $e=aU+bT-
\sum_i c_i E_i$ is a class.

\noindent {\it Properties 1 and 2.} We will prove the first two
properties together. 
As we have mentioned, the -Id automorphism is realized by an
orientation-preserving diffeomorphism.
Therefore we can assume that $a\geq 0$. 

The easier case is when $a=0$. In this case, 
 $ e\cdot e=-\sum_i  c_i^2$. If  $e$ has non-negative square, 
then $c_i=0$ for each $i$ and 
$e$ is simply the  class 
$bT$, which is certainly reduced.
 If  $e$ has square $-1$. 
then  
$c_i=\pm 1$ for some $i$ and $c_j=0$ for any $j\ne i$.   
Consider the reflection along $[b/2]T-E_i$ which maps 
$e$ to $E_i$ or $T+E_i$. When $n\geq 1$,  the reflection along
$T+E_1-E_2$
maps $T+E_i$ to  the class $E_2$. 
If  $e$ has square $-2$, we have 
$c_i=c_j=1$ for some $i\ne j$ and $c_k=0$ for any $k$ different from
$i$ and $j$.

When
$a$ is strictly positive,   
we will show that under an 
orientation-preserving diffeomorphism
which is a composition of reflections along classes represented by
embedded spheres with square $-1$, $e$ 
is equivalent to a class $\tilde e=aU+\tilde b-\sum_i \tilde c_i E_i$
with $a\geq \tilde c_i$ for each $i$.  
For any  $r_i\geq 0$ and $\epsilon_i=\pm1$ to be determined, 
it is easy to see, via the tube construction,  that $\mu_i=r_iT+\epsilon_i E_i$
is represented by an embedded sphere with square $-1$.   
Therefore, the  reflection along $\mu_i$ is realized by an
orientation-preserving diffeomorphism. 
Since $e\cdot \mu_i=ar_i+\epsilon_i c_i$, under the reflection, 
$$c_i\longrightarrow c_i'=c_i-2ar_i\epsilon_i-2c_i=-c_i-2ar_i\epsilon_i$$ 
 and $a$ is invariant. We first assume that $a$ is positive. In order for $|c_i'|\leq a_i$, we find that $r_i$ and $\epsilon_i$ should satisfy
$$-{1\over 2}-{c\over 2a}\leq r_i\epsilon_i\leq {1\over 2}-{c\over 2a}.$$
Clearly, such a pair ($r_i,\epsilon_i$)
exists, and there is a unique solution when $c/a$ is not an odd  integer,
and there are two solutions when $c/(2a)$ is an odd integer.
By applying this process for each $i$, we obtain a desired class
$\tilde e$. 
Notice that
$\tilde e$ is equivalent to a reduced class under trivial
automorphisms. So we have proved that  $e$ is equivalent to
a reduced class if $a>0$.

\noindent{\it Property 3.}  It is a immediate consequence of the fact (see [Bi] or [LLiu1])
that $${\cal E}_{K_0}=\{E_1, ..., E_n, T-E_1, ..., T-E_n\}.$$
Indeed, $e\cdot E_i=c_i$ and $e\cdot (U-E_i)=a-c_i$,  
both of which are  positive because $e$ is reduced. 

\noindent{\it Property 4.} Consider  symplectic forms 
$\omega_{\epsilon}=U+T-\sum \epsilon E_i$. Their  
 canonical
class 
is $K_0=-2U+(2g-2)T+\sum_i E_i$, and for $\epsilon$ small,
$\omega_{\epsilon}\cdot e >0$ for any reduced class $e=aU+bT-\sum_i c_i E_i$.
Therefore, $K_0$ is in ${\cal K}_e$. 

\noindent{\it Property 5.}
Suppose now $e=aU+bT-\sum_i c_i E_i$ is a reduced class
with non-negative square.
Let $\tau$ be the class of a symplectic form which
is positive on $e$. Since both $a$ and $b$ are non-negative and one of
them is positive, by the light cone lemma,
$\tau$ must have a positive $U$ term as well.    
Therefore any symplectic canonical class in ${\cal K}_e$ is of the form
$K=-2U+bT+\sum s_i E_i$ with $s_i$ being odd. Since $K^2=8-8g-n$, 
$b=2g-2-[(\sum_i c_i^2-n)/4]$. 
Thus we have $K_0-K=[(\sum_i s_i^2-n)/4]T +(1-s_i)E_i$, and consequently
$$\eqalign{(K_0-K)\cdot e&=  a\sum_i (s_i^2-1)/4 + \sum_i (1-s_i)\cr
&=\sum_i (1-s_i)[{(3-s_i)+(1-a)(1+s_i)\over 4}].\cr}$$
Let $S_i=(1-s_i)[{(3-s_i)+(1-a)(1+s_i)}]$.
We will show that $S_i\geq 0$ for each $i$. 
 When $s_i\geq 3$, the two factors of $S_i$ are both non-positive, so
$S_i$ is non-negative. When $s_i=1$, $S_i=0$. Finally, when 
$s_i\leq -1$,  the two factors are both non-negative,
and therefore $S_i$ is non-negative. We have finished the proof
of property 5 for a non-minimal ruled manifold and
hence the proof of Lemma 3.4.

 We will now prove a rather general result relating the
symplectic genus and  the minimal genus of a reduced class, 
using Taubes' equivalence between Seiberg-Witten invariants
and Gromov-Taubes invariants ([T2]). 
Let us first provide some background of this 
equivalence (see e.g. [LLiu1] and [T2]). 

Recall that Seiberg-Witten invariants are defined on Spin$^c$ structures.
For manifolds without torsion-free homology group, like rational
and irrational ruled manifolds, the Spin$^c$structures correspond to
characteristic classes. For this reason, we will simply
speak of the Seiberg-Witten invariants of the characteristic
classes. Suppose $K$ is a symplectic canonical class,
 then any class 
of the form $-K+2e$ is a characteristic class.  The Seiberg-Witten
invariant of $-K+2e$ is defined when its Seiberg-Witten
dimension $-K\cdot e+e\cdot e$ is non-negative. For
manifolds with $b^+=1$, the Seiberg-Witten invariants also
depend on the chambers. In the presence of a symplectic form $\omega$, 
there is an $\omega-$symplectic chamber. On such a manifold,
the Gromov-Taubes invariant of a class $e$ is a suitable
count of  $\omega-$symplectic
surfaces representing $e$. The surface is not required to be
connected, but is required to be embedded and any component
with negative square is a $\omega-$symplectic sphere with
square $-1$.   

 When $K$ is the symplectic canonical
class of $\omega$, a fundamental theorem of Taubes
states that, if the Seiberg-Witten invariant of $-K+2e$ in the 
$\omega-$symplectic chamber is nontrivial, then (i)
$e$ is represented by a $J-$holomorphic
curve (possibly singular) for any $\omega-$compatible
almost complex structure $J$; (ii)  the Seiberg-Witten invariant is 
the same as the Gromov-Taubes invariant of $e$ provided
$e\cdot E\geq 0$ for any $E\in {\cal E}_K$. 

\noindent{\bf Proposition 3.5}. Let $M$
be a non-minimal rational or irrational ruled manifold with a standard
decomposition and a standard basis.  Suppose $e$ is a reduced
class. 
If $e\cdot e$ is no less than
$\eta(e)-1$, then $e\cdot e\geq 0$ and $e$ is represented
by a symplectic surface. Moreover,
 if $e$  is either a class of positive square or
a primitive class with square 0,   
$e$ is represented
by a connected   symplectic surface, and therefore its minimal genus
 is given by  its symplectic genus.

\noindent{\it Proof}. We will first prove  that $e$ is represented by
a symplectic surface.
 By the definition of
the symplectic genus and Lemma 3.4(4)
$$K_0\cdot e+e\cdot e\leq 2\eta(e)-2.$$
Therefore, under our assumption, the Seiberg-Witten dimension 
of the class $-K_0+2e$ satisfies  
$$-K_0\cdot e+e\cdot e\geq 2(e\cdot e +1-\eta(e))\geq 0.$$ 
Now we divide the proof into two cases. 

In the case of rational manifold, 
for a symplectic from $\omega$ with $-K_0=3H-\sum_i E_i$ as its  canonical 
class, it is shown in [LLiu2] that 
$H$ is represented by an embedded $J-$holomorphic sphere
for a generic almost complex structure $J$ compatible with $\omega$. 
Since the  reduced class $e=aH-\sum b_i E_i$ has a positive $a$ term,
 $(K_0-e)$ has a negative $a$ term and so $(K_0-e)\cdot H<0$.
Therefore, $K_0-e$ is not represented by a $J-$holomorphic curve,
  because  the intersection number
of two distinct $J-$holomorphic curves is non-negative. 
So the Seiberg-Witten invariant of $-K_0+2(K_0-e)=K_0-2e$ is trivial by
the result of Taubes. By the symmetry of Seiberg-Witten
invariants (see Lemma 2.3 in [LLiu1]), the Seiberg-Witten invariant
of $-K_0+2e$ in the non-symplectic chamber is trivial. 
By the wall crossing formula of Seiberg-Witten invariants
(see [KM] and Lemma 3.3 in [LLiu1]),
the Seiberg-Witten
invariant $SW_{\omega}(-K_0+2e)$  in the $\omega-$symplectic chamber 
is non-trivial.
Since $e$ is reduced, by Lemma 2.3(3),
we have $e\cdot E\geq 0$ for any $E\in {\cal E}_{K_0}$.
Thus, $e$ is represented by an embedded symplectic surface 
by the result of Taubes.

In the case of irrational ruled manifold,
by [LLiu2], for any symplectic form $\omega$ with
$K_0$ as its canonical class, $T$ is represented 
by a $J-$holomorphic sphere for a generic $\omega-$compatible
almost complex structure $J$.
Since a reduced class  has a positive $U$ term and $U\cdot T=1$,
 we can show that $K_0-e$ has trivial Seiberg-Witten invariant
in the $\omega-$symplectic chamber.
Applying Lemma 2.3 and Lemma 3.3 in [LLiu1] as above,
and notice that $-K_0+2e$ has a positive
$U$ term and that the class $\gamma$ in Lemma 3.3 in [LLiu1] is just the
class $T$ here, we find that the Seiberg-Witten invariant of
$-K_0+2e$ is nontrivial.   
Taubes's result and Lemma 2.3(3) then can be applied to show that $e$ 
is represented by an embedded symplectic surface. 

We have shown that $e$ is represented by a symplectic surface.
This surface may have many components. Any component
with negative square is a symplectic sphere with square $-1$.
However, since $e\cdot E\geq 0$ for any $E\in {\cal E}_{K_0}$,
no such component exists. 
Thus, $e$ is represented by a symplectic surface
the components of which all have non-negative square,
and therefore $e\cdot e$ is non-negative.  
If $e\cdot e>0$, there can  
only be one component by the light cone lemma. 
If $e\cdot e=0$,  again by the light cone lemma,
there might be several components, all of which  are  multiples of  the
same class. All the multiplicities 
have to be  one because of the adjunction formula. Thus, if $e$ is primitive,
there is only one component. 
The proof is complete.


Notice that, as an immediate consequence of   Proposition 3.5, 
  the symplectic genus of certain reduced class 
is non-negative. In fact, this weaker assertion holds in 
much greater generality.

\noindent {\bf Lemma 3.6}. Let $M$ be a non-minimal rational or 
irrational ruled manifold with a standard decomposition
and a standard basis.

\noindent 1. The symplectic genus of any class with positive square
or a primitive class
with square 0  is non-negative.

\noindent 2.  Any class with square $-1$ or $-2$ has non-negative 
symplectic genus. In addition, the  classes which are equivalent
to reduced classes
 have positive symplectic genus, and those which are not equivalent
to  reduced classes
have  symplectic genus 0.
 

\noindent{\it Proof}.  
Let $e$ be a class with square at least $0$ and
equivalent to a reduced class $e'$. 
Due to Lemma 3.2(1), $e$ and $e'$ have the same symplectic genus.
Suppose that the symplectic genus of $e$ is negative, then $e\cdot e\geq -1\geq
\eta(e)-1$.   
By Proposition 3.5, $e'$ is represented by a 
symplectic surface and hence the connected symplectic genus is non-negative. 
This is a contradiction.

When $e\cdot e=-1$, by Lemma 3.4(2), $e$ is either equivalent
to a reduced class, or  equivalent  to 
$E_1, H-E_1-E_2$ or $T-E_1$. 
It is easy to see that
 $E_1, H-E_1-E_2$ or $T-E_1$ are all spherically representable and
have symplectic genus zero. 
Suppose $e$ is a reduced class  and
$\eta(e)\leq 0$. Since $e\cdot e=-1$,
 it satisfies the assumption of Proposition 3.5,
 and  we can conclude that $e\cdot e\geq 0$. This  contradicts with
our assumption. Therefore, by Lemma 3.2(1),
 any class equivalent to a reduced class  has positive symplectic genus. 

For the case of a class of square $-2$, the same argument as in
the previous paragraph proves that
the symplectic genus can not be smaller than zero. What we still need
to show is that there does not exist
 any  reduced class $e$ with symplectic genus 0.
Suppose $e$ is such a class. Then by definition there is a
$K\in {\cal K}$ such that $K\cdot e=0$.
In light of Lemma 3.4(4), it is also  necessary that $K_0\cdot e\leq 0$.
 We first exclude the case $K_0\cdot e <0$.  
 Let $K'$ be a symplectic canonical class
such that ${\cal C}_{K'}$ and ${\cal C}_{K_0}$ are
in the same component of the positive cone ${\cal P}$. 
Notice that the argument in
the proof of Lemma 3.4(5) actually proves that 
$K'\cdot e\leq K_0\cdot e$ for any reduced class $e$. 
Therefore all such $K'$ satisfies $K'\cdot e<0$.  It is clear that
any symplectic canonical
class is either a $K'$ or a $-K'$.  Thus, there is no
symplectic canonical class $K$ satisfying $K\cdot e=0$. 
This contradiction leaves the case $K_0\cdot e=0$ as the
only possibility. In this case, the reflection $R^{e}$
along $e$ preserves ${\cal E}_{K_0}$ since it preserves
$K_0$. So, if $F\in {\cal E}_{K_0}$, then 
$F'=R^{e}(F)\in {\cal E}_{K_0}$, and 
$F'=F+(e\cdot F)e.$ By [LLiu2], for any symplectic 
form $\omega\in {\cal C}_{K_0}$, $F$ and $F'$
are both represented by smooth $J-$holomorphic spheres
for some generic $\omega-$compatible almost complex structure
$J$, we have $F\cdot F'\geq 0$ by the positivity of
intersection. This fact, together with Lemma 3.4(3),
leads to the following contradiction
$$-1=F'\cdot F'=F\cdot F' + (e\cdot F) e\cdot F'\geq 0.$$ 
The lemma is proved. 

We are ready to prove Theorem B. In fact, we will prove
the following  more general result.  

\noindent{\bf Theorem B'}. Let $M$ be a rational or irrational ruled four$-$manifold.  
Suppose $e$ is a class with square at least $-1$, and in the case
that $e$ has square one, we further assume that 
$e$ is a primitive class. Then 
its symplectic genus is non-negative and there is an algorithm to
calculate its symplectic genus.
Furthermore, if $e\cdot e\geq \eta(e)-1$, then $e$ is represented
by a connected symplectic surface, and therefore its minimal 
genus coincides with its symplectic genus.  

\noindent{\it Proof}. When $M$ is minimal,
$M$ is either $CP^2, S^2\times S^2$ or an $S^2-$bundle over a 
Riemann surface. The minimal genus problem for these manifolds has been
completely solved in [LiL3-4].  

When $M$ is non-minimal, with a choice of  a standard decomposition and
a standard basis, we can define reduced classes. 
Suppose $e$ is a class satisfying the conditions of Theorem B'.
By Lemma 3.4, under an algorithm, $e$ can be transformed to a 
reduced class $\tilde e$
or a class $e'$ which can be represented by a symplectic sphere.
 
Since $e\cdot e=\tilde e\cdot \tilde e$, and $\eta(\tilde e)=\eta(e)$
by Lemma 3.2(4), we see that $\tilde e$ satisfies the conditions of
Proposition 3.5 and Lemma 3.6(1).   By Lemma 3.4(4),
the symplectic genus of $\tilde e$ can be calculated with 
a simple formula. And by Lemma 3.6(1), the symplectic genus
of $\tilde e$ is non-negative. Finally, by Proposition 3.5, $\tilde e$
is represented by a connected symplectic surface.
The proof of Theorem B' is complete.   
 
\noindent{\bf Theorem A}. Let $M$ be a smooth, closed  oriented
$4-$manifold with non-empty symplectic cone and 
$b^+(M)=1$. Then the symplectic genus of any class of positive square is
non-negative, and it coincides with the minimal
genus for  any sufficiently large multiple of such a class.

\noindent{\it Proof}. 
   In the rational and irrational
ruled cases, by Theorem B',
every class with positive square has positive symplectic genus.
If $M$ is neither 
 rational nor irrational ruled,  we examine the  minimal case first. 
Given a class $e$ with positive square and
 a symplectic form $\omega$, by the light cone lemma, either
$\omega\cdot e>0$ or $-\omega\cdot e>0$. Let us assume that we are
in the first situation.  
By a result in [Liu],  $K(\omega)\cdot \omega\geq
0$ if $K$ is the canonical class of a symplectic form $\omega$. 
Then, by the light cone lemma,  $K(\omega)\cdot e\geq 0$.  
Thus it follows directly from inequality (3.2) that 
the symplectic genus of $e$ is positive. 
For the non-minimal case, we claim that 
one can find $K\in {\cal K}_e$ such that $K\cdot e\geq 0$ if $e\cdot e\geq 0$. 
The positivity of $\eta(e)$ then follows immediately from it.
Suppose $M=N\# n{\overline CP}^2$ is the (unique) minimal
reduction of $M$. Let $E_1, ..., E_n$ be the generators of
$H^2$ of the $n$ $\overline {CP}^2$.  
Write $e=e_m-r_1E_1-...-r_n E_n$, where $e_m$ is the pull back of a class
in $H^2(N;{\bf Z})$ also denoted by $e_m$. Pick a symplectic form 
$\omega_m$ on $N$ such that $\omega_m\cdot e_m>0$.   
Let $K_m$ be a symplectic canonical
class of $\omega_m$. Then, as above, we have $e_m\cdot K_m\geq 0$. 
By the blow up construction, for sufficiently small $\epsilon$, 
the class $[\omega_m]-\epsilon E_1-...-\epsilon E_n$
is realized by a symplectic form on $M$ with symplectic canonical class 
$K_m+E_1+...+E_n$. Applying the reflections along the $E_i$, we see
that $[\omega_m]\pm \epsilon E_1\pm ...\pm \epsilon E_n$ are
realized by symplectic forms with symplectic canonical classes
$K_m\mp E_1\mp..., \mp E_n.$ 
For possibly smaller $\epsilon$, the pairing between
$e$ and $[\omega_m]\pm \epsilon E_1\pm ...\pm \epsilon E_n$
is positive.     
Therefore,  any symplectic canonical class of the
form $K=K_m \pm E_1\pm...\pm E_n$ is in ${\cal K}_e$. 
Since $K_m\cdot e_m\geq 0$, by choosing $E_i$ or $-E_i$
appropriately, we can easily find a $K\in {\cal K}_e$ such that
$K\cdot e\geq 0$. 

The last statement of the theorem (for a class $e$ satisfying
$e\cdot E\ne 0$ for any $E\in {\cal E}$) 
 is a direct consequence of the following two results, together with
Lemma 3.2(1). 
One  result is in [LLiu1] that 
$$C_M=\{e\in {\cal P}|0<|e\cdot E|\hbox { for all $E\in {\cal E}_M$ }\}.$$
The other is due to Donaldson (see [D]). It states that,
 for any sufficiently large integer $N$,  $N[\omega]$
can be represented by connected symplectic submanifolds.
Now suppose that $e\cdot E=0$ for some $E\in {\cal E}$. 
By the result in [L1], there exists a symplectic form $\omega$ such that
$E$ is represented by an $\omega-$symplectic sphere. Blowing down that sphere,
we obtain a new symplectic manifold $M'$. There is a class $e'$ in $M'$
which is pulled back to $e$. It is easy to see that 
$m(le')\geq m(le)$ and $\eta(le')=\eta(le)$ for any integer $l$.
If $e'\cdot E'\ne 0$ for any $E'\in {\cal E}_{M'}$, then $\eta(le')=m(le')$  
for sufficiently large $l$. Therefore $\eta(le)=\eta(le')=m(le')\geq
m(le)$. Together with Lemma 3.2(1) we arrive at the conclusion
that $\eta(le)=m(le)$. If there is still a class $E'\in {\cal E}_{M'}$
such that $e'\cdot E'=0$, we can continue the process above. However, this
process can only be repeated finitely many times. 
The proof of Theorem A is complete. 

We remark that, using some of the arguments in
[LLiu1],  in fact we are able to  get an effective estimate  on
how large a multiple $N$ is allowed in the last statement of Theorem A.   
Here we just mention, in the case of a minimal manifold with $b^+=1$
which is
neither rational nor irrational ruled, it suffices to take 
$N=2|e\cdot K|/e^2$, where $\pm K$ are the only
two symplectic canonical classes.  
In particular,  when a manifold with $b^+=1$ has a torsion symplectic
canonical class,
we are able to conclude that the minimal genus
of  every class $e$ with positive square coincides with its symplectic genus
(which is simply $(e\cdot e)/2+1$). Such manifolds include the 
Enriques surface,
hyperelliptic surfaces, any torus bundle over torus which has
$b^+=1$. In addition, from the results in [LiL4], [Li1] and [Kr2],  
manifolds with the property that two genera coincide for
any class of positive square include minimal irrational ruled manifold, 
rational manifold with $b^-\leq 9$ and the product of
a circle with a fibered $3-$manifold $Y$ with $b_1(Y)=1$.

We close this section with another remark.
There are classes of positive square, which do 
not satisfy the conditions of Theorem
B but still have the same symplectic genus and minimal genus.
Some of them are actually   
represented by connected symplectic surfaces. For any positive integer
$a$ bigger than 4,
consider the reduced class $aH-\sum_{i=1}^{a^2-1} E_i$. Its 
 square is 1 and symplectic genus  $(a^2-3a)/2$. If we blow up
$a^2-1$ points on a smooth curve of degree $a$, then
the proper transform is a smooth curve in this given class. 
Others, including some classes in the non-trivial 
$S^2-$bundles over Riemann surfaces
  are  not known to be represented by connected
symplectic surfaces. 
To deal with such classes, we may need to 
find more  constructive techniques as in [LiL3-4] and [Li1].

\bigskip

\noindent{\bf \S4. The classes represented by spheres}
\medskip
In this section we determine  the set of classes 
represented by spheres and the orbits of Diff$(M)$ on this set. 
We start with 

\noindent{\bf Theorem C}. Let  $M$ be a  rational or irrational ruled manifold and
$e\in H^2(M)$ be a class with square at least $-1$. 
  If
$\eta(e)=0$, then  
PD$(e)$ is represented by a smoothly embedded sphere.
Furthermore, if
PD$(e)$ is represented by a smoothly embedded sphere,
then either $\eta(e)=0$ or $e$ is a non-primitive class of square
zero with $e=pe'$ and $\eta(e')=0$.

\noindent {\it Proof}. 
Suppose a standard decomposition and a standard basis are
given. 
Let us first deal with the case  $e\cdot e\geq 0$. 
Now suppose $m(e)=0$. Then by Lemma 3.2(1), $\eta(e)\leq 0$.
By Lemma 3.6(1), $\eta(e)=0$ unless 
$e$ is a divisible class with square zero, ie. $e=pe'$ for some $p\geq 2$
and some $e'$ with $e'\cdot e'=0$. Since $\eta(pe')=\eta(e)\leq 0$,
 by Lemma 3.2(3), 
$\eta( e')$ 
can not be positive.  
In this case, by Lemmas 3.2(4) and 3.4(1), there is a reduced primitive
class 
$\tilde e'$ with the same square,  the same symplectic genus and
the same symplectic minimal genus as $e'$.    
Since $\tilde e'$ is primitive and reduced with   $\tilde e'\cdot \tilde
e'= 0>\eta(\tilde e')-1$,  we can apply Proposition 3.5 to conclude 
that 
$\eta(\tilde e')$ coincides with $m(\tilde e')$. Since 
$\eta(\tilde e')\leq 0$, and $m(\tilde e')\geq 0$ by definition,
both of them are equal to zero.  
 Therefore, in this case,   $e'$  has symplectic genus zero as well. 
 
Suppose the symplectic genus $\eta(e)$ is zero.
Again, there is a reduced class 
$\tilde e$ with the same square, the same divisibility,
the same symplectic genus and
the same symplectic minimal genus. Applying Proposition 3.5 to $\tilde e$,  
 together with
 Lemma 3.2(3), which excludes the case when $\tilde e$ is a divisible
class with square zero, 
we conclude that 
$m(\tilde e)=0$. Therefore, $m(e)$ is zero as well. 

Finally we deal with the case
that  $e\cdot e=-1$.
By Lemma 3.6(2),  either 
$e$ has positive symplectic genus, or
$\eta(e)=0$ and $e$ is spherically representable.
When $\eta(e)>0$, $e$ is not spherically representable by a sphere 
due to Lemma 3.2(1). Thus, $e$ is spherically representable
if and only if $\eta(e)=0$.
The proof is finished.

For the convenience of the proof of Theorem D, we
state the following corollary.

\noindent {\bf Corollary 4.1}. Let $M$ be a rational or irrational
ruled $4-$manifold.   Suppose $e$ is a class with positive square   or
a primitive class with square zero,  
the following statements are equivalent:

\noindent
1. $e$ is represented by a smoothly embedded sphere.

\noindent 2. $\eta(e)$ is zero.  

\noindent 3. $e$ is represented by a symplectically embedded sphere 
with respect to some symplectic form.

\noindent{\it Proof}. The equivalence of the first two statements
follows from Theorem C.  The equivalence of last two
statements follows directly  from Proposition 3.5 and Lemma 3.4(1).

We remark  that Corollary 4.1  holds for classes with 
square $-1$ and $-2$ as well.
 
Having determined the set ${\cal SPH}_{\geq -1}(M)$, we are going to 
classify the orbits
of Diff$(M)$ on this set.  
We begin with the difficult case when $M$ is rational.

\noindent {\bf Theorem 4.2}. Let $M$ be a rational manifold
with a standard decomposition and a standard basis.
Then the following classes are spherically representable: 

\noindent 1. $2H$,

\noindent 2. $ (k+1)H-kE_1, \,\, k\geq 0,$

\noindent 3. $(k+1)H-kE_1-E_2, \,\, k\geq 1,$

\noindent 4. $kH-kE_1, \,\, k\geq 1.$

\noindent Moreover, up to  diffeomorphisms, 
any spherically representable class
with non-negative squares is equivalent to one of the above. 

\noindent{\it Proof}.  The first claim is well known. We just give a sketch here.
$H$ and the $E_i$ are all spherically representable. Moreover,
the spheres representing them can be chosen to be pairwise disjoint. The first claim now  follows from the elementary
fact: if  $A_1$ and $A_2$ are represented by two spheres
which  intersect at most at one point, then  $A_1+A_2$ is spherically
representable. 

To prove the last claim, we need the following two results.

\noindent{\bf Proposition 4.3}. Up to automorphisms of $H^2$, 
the set of spherically representable classes
with non-negative square are given as above.

\noindent {\bf Lemma 4.4}. Let $\omega$ be a symplectic form
with symplectic canonical class $K$. 

\noindent 1. Any class $R$ with positive square 
and represented by an $\omega-$symplectic sphere is
in $\hat {\cal C}_K$ and satisfies $R\cdot(- K)>0$.

\noindent 2. Any  $R$ with square 0 and represented by
an $\omega-$symplectic sphere is in $\overline {{\cal C}_K}\cap {\cal B}$
and satisfies $R\cdot (-K)>0$.

\noindent{\it Proof of Proposition 4.3}. The classes in the first three
cases have positive square.  
In this case, the claim was proved by Kikuchi 
in [K1]. 
For classes of square 0, this was implicitly shown in [Li2].  
In fact, it was shown in [Li2] that 
 if $y$ is  a primitive class represented by an embedded sphere, then there exists $x, z_1, \dots, z_{n-1}$ such that
$$H^2(M;{\bf Z})=\pmatrix{0&1\cr 1&t\cr}\oplus (n-1)(-1)$$
with respect to the basis $<y, x, z_1,\dots,  z_{n-1}>$.
If $t$ is odd, let $\tilde x=x-[(t-1)/2] y$, 
$\tilde z_i=z_i$ for $i=1,...,n-1$;
if $t$ is even, let $\tilde x=x-[(t-2)/2]y+z_1$, $\tilde z_1=z_1+y,$ 
$\tilde z_i=z_i$ for $i=2,..., n-1.$
Then, with respect to the new basis $<y, \tilde x, \tilde z_1,...,
\tilde z_{n-1}>, $
 $$H^2(M;{\bf Z})=\pmatrix{0&1\cr 1&1\cr}\oplus (n-1)(-1).$$
Since $H^2(M;{\bf Z})$ has the same decomposition with respect
to the basis $<H-E_1, H, E_2, ..., E_n>$,
 there is an automorphism of $H^2(M;{\bf Z})$ sending $y$ to 
$H-E_1$. The non-primitive case follows immediately.

\noindent {\it Proof of Lemma 4.4}. For any $\omega-$compatible almost complex structure
$J$, the $\omega-$symplectic sphere representing $R$ can be taken $J-$holomorphic.
Moreover, 
for a generic $\omega-$compatible almost complex structure
$J$, any $E\cdot {\cal E}_K$ is represented by a smooth
$J-$holomorphic sphere.  Then   
$R\cdot E\geq 0$ for any $E\in {\cal E}_K$ by the positivity
of intersection of pseudo-holomorphic curves. Thus, when $R$ has positive
square, it  is  $\hat {\cal C}_K$, and  when $R$ has square zero, it is in the
$\overline {{\cal C}_K}\cap {\cal B}$.
In either case, by the adjunction formula,
$R\cdot (-K)=2+R\cdot R\geq 2$. 
  The lemma is proved. 

We now continue  the proof of Theorem 4.2. 
Suppose  $R$ is a class with non-negative square
and represented by a sphere. Let $R'$ be a class in the list of Theorem 4.2
 with the same square and the same divisibility.
By Proposition 4.3, there exists an  automorphism
$\phi$ such that  $\phi^*(R)=R'$.
By Corollary 4.1, both $R$ and $R'$ are represented by 
symplectic spheres (possibly with respect to different symplectic forms). By Lemma 4.4 and  Theorem 2.8, we see that $\phi$ is realized by 
an orientation-preserving
 diffeomorphism. The proof of Theorem 4.2 is complete.
 
We remark that,  for a rational manifold with $b^-\leq 9$,
Theorem 4.2 has been proved in [K] and [Li1-2]. In fact, in this case,
it follows immediately from the fact that  
every automorphism is realized by an orientation-preserving
diffeomorphism.  
Kikuchi ([K2]) also  conjectured that Theorem 4.2 would hold in general.

Now we are ready to prove Theorem D. 

\noindent{\it Proof of Theorem D}. When $M$ is neither rational
nor irrational ruled, $M$ has a unique minimal reduction
$M=N\# n\overline{CP}^2.$ 
Let $E_1,..., E_n$ be the generators of the $\overline 
{CP}^2$. The only spherically representable classes
with square at least $-1$ are $\pm E_1, ..., \pm E_n$.
They are carried to each other by trivial automorphisms.

On a minimal  irrational ruled manifold $M$, among all classes
with square at least $-1$, up to sign, 
there is a unique primitive 
class which  is spherically representable (see [LiL4]). Since $-$Id is
in D$(M)$, again there is a unique orbit when the square and
the divisibility are fixed. 

A minimal rational manifold
is either $CP^2$ or $S^2\times S^2$. Let $H$ be a generator of
$H^2(CP^2;{\bf Z})$.  $\pm H$ and $\pm 2H$ are the only spherically representable classes,
with square $1$ and $4$ respectively (see [KM]). Since -Id$\in$ D$(CP^2)$,
there is only one orbit when the square is fixed. 
Let $x$ be the class represented by $S^2\times \{pt\}$
and $y$ be the  class represented by $\{pt\}\times S^2$.  
For each even number $2l$, there are  four spherically representable classes
with square $2l$: $\pm (x+ly)$ and $\pm (lx+y)$ (see [Ku]). 
 Since $-$Id and 
the automorphism switching the two factors are in D$(M)$,
the uniqueness of the orbits for fixed square is  obvious.

Finally, let us consider the non-minimal manifolds.
In the rational case, suppose 
$M$ is given a minimal reduction $M=CP^2\#n\overline {CP}^2$
and a standard basis; and in the irrational case, suppose
$M$ is given  a minimal reduction of the form
$S^2\times \Sigma_h\# (n-1)\overline {CP}^2$ and a standard basis.
We first treat the case of negative square.

\noindent{\bf Proposition 4.5}. Diff$(M)$ has one orbit on ${\cal SPH}_{-1}(M)$
 when $n=1$ and $n\geq 3$, in the exceptional case $n=2$,
Diff$(M)$ has two orbits, one ordinary and one characteristic. 
Diff$(M)$ has one orbit
on ${\cal SPH}_{-2}(M)$ when $n=2$ and $n\geq 4$, 
in the exceptional case $n=3$,
Diff$(M)$ has two orbits, one ordinary and one characteristic.

\noindent{\it Proof}. This follows from Lemma 3.4(2) and Lemma 3.6(2).  

For an irrational ruled manifold, the only spherically representable classes with non-negative square
are $\pm kT$.

From Theorem C we can  list all the possible orbits with non-negative square.
For a given square, when there are more than one orbit, they are distinguished
by divisibility. 

\noindent $s=0$: Diff$(M)$ has 
infinitely many  orbits on ${\cal SPH}_s(M)$, represented by $k(H-E_1)$,

\noindent $s\geq 1$ and odd: Diff$(M)$ has 
one orbit on ${\cal SPH}_s(M)$, represented by $[(s+1)/2]H-[(s-1)/2]E_1$,

\noindent  $s=2$ or $s\geq 6$  and even: Diff$(M)$ has one orbit on ${\cal
SPH}_s(M)$ if  $l\geq 2$, represented by $[(s+2)/2]H-([s/2)]E_1-E_2$,

\noindent $s=4$: Diff$(M)$ has one orbit on ${\cal SPH}_s(M)$ if $l\leq 1$,
represented by $2H$; 
two orbits if  $l\geq 2$, represented by $2H$ and $3H-2E_1-E_2$.

Theorem D is thus proved for all cases.



\bigskip
\noindent {\bf Reference}
\medskip

\noindent [Bi] P. Biran, Symplectic packings in dimension 4, 
Geom. and Funct. Anal. 7 (1997), no.3. 420-437. 


\noindent [D] S. Donaldson, Symplectic submanifolds and almost-complex 
geometry, J. Differential Geom. 44 (1996), no.4. 666-705.

\noindent [FM1] R. Friedman and J. Morgan, On the diffeomorphism types 
of certain algebraic surfaces, 
J. Differential Geom. 27 (1988), no.3 371-398.

\noindent [FM2] R. Friedman and J. Morgan, 
Algebraic surfaces and Seiberg-Witten invariants. J. Algebraic Geom. 6
(1997), no. 3, 445--479.

\noindent [FS] R. Fintushel and R. Stern, Immersed spheres in $4-$manifolds
and the immersed Thom conjecture, Turkish J. Math. 19(2) (1995), 145-157.

\noindent [G] H. Z. Gao, Representing homology classes of $4-$manifolds,
Topology and its Application, 52(2) (1993), 109-120.  


\noindent [K1] K. Kikuchi, Positive 2-spheres in 4-manifolds of signature
$(1, n)$, Pacific J. Math., 160 (1993), 245-258. 

\noindent [K2] K. Kikuchi, Personal communication.

\noindent [KM] P. Kronheimer and T. Mrowka, The genus of embedded
surfaces in the projective plane, Math. Res. Letters 1 (1994), 797-808.

\noindent [Kr1] P. Kronheimer, Embedded surfaces and gauge theory in 
three and four dimensions, Survey in differential geometry,
Vol. III (Cambridge, MA, 1996), 243-298, Int. Press, Boston,
MA, 1998.

\noindent [Kr2] P. Kronheimer, Minimal genus in $S^1\times M^3$,
Invent. Math. 135 (1999), no.1, 45-61.

\noindent [Ku] K. Kuga, Representing homology classes of $S^2\times S^2$,
Topology 23 (1984), 133-137.

\noindent{[La1]} T. Lawson, { Smooth embeddings of 2-spheres in
 4-manifolds}, Expo. Math. 10 (1992), 289-309;

\noindent{[La2]} T. Lawson, { The minimal Genus problem}, Expo. Math. 15(1997), 385-431;

\noindent [L1] T. J. Li, Smoothly embedded spheres in symplectic
four manifolds,  Proc. Amer. Math. Soc. 127 (1999), no. 2, 609--613.

\noindent [L2] T. J. Li, Fiber sums of Lefschetz fibrations, in preparation.

\noindent [Li1] B. H. Li, 
Representing nonnegative homology classes of $CP^2\#{\overline {CP}^2}$
by minimal genus embeddings, Tran. Amer. Math. Soc.
352 (2000), no. 9. 4155-4169.

\noindent [Li2] B. H. Li, Embeddings of surfaces in $4-$manifolds (II),
Chinese Science Bull.
36 (1991), 2030-2033.  

\noindent [LLiu1] T. J. Li and A. K. Liu, Uniqueness of symplectic canonical class, surface cone and symplectic cone of $4-$manifolds with $b^+=1$,
preprint.

\noindent [LLiu2] T. J. Li and A. K. Liu, Symplectic
structures on ruled surfaces and a generalized adjunction inequality, Math.
Res. Letters 2 (1995), 453-471.


\noindent [LiL2] B. H. Li and T. J. Li, Smooth minimal genera for
small negative classes in $CP^2\#n{\overline{CP}^2}$ when $n\leq 9$, preprint.

\noindent [LiL3]  B. H. Li, T. J. Li, {Minimal genus smooth embeddings in
$S^2\times S^2$ and $CP^2\# n\overline{CP}^2$ with $n\leq 8$},
Topology, 37 (1998), 573-594.

\noindent [LiL4] B. H. Li, T. J. Li, { Minimal genus embeddings of
surfaces in $S^2$-bundles over surfaces}, Math. Res. Lett. 4(1997), 379-394.

\noindent [Mc1] D. McDuff, The structure of rational and ruled
symplectic
$4$-manifolds. J. Amer. Math. Soc. 3 (1990), no. 3, 679--712.

\noindent [Mc2] D. McDuff, Lectures on
Gromov invariants for symplectic 4-manifolds, 
With notes by Wladyslav Lorek. NATO Adv. Sci.
Inst. Ser. C Math. Phys. Sci., 488, Gauge theory and symplectic geometry
(Montreal, PQ, 1995), 175--210, Kluwer Acad. Publ., Dordrecht, 1997.

\noindent [Mc3] D. McDuff, Immersed spheres in symplectic $4-$manifolds, 
Ann. Inst. Fourier (Grenoble) 42 (1992), no.1-2, 369-392.


\noindent [MS] D. McDuff and D. Salamon, A survey of symplectic
$4-$manifolds with $b^+=1$, Proc. Gokovo conference, 1995, 47-60. 

\noindent{[OS]} P. Oszvath, Z. Szabo, The symplectic Thom conjecture,
Ann. of Math. (2) 151 (2000), no. 1, 93-124.


\noindent [T1] C.H. Taubes, {The Seiberg-Witten invariants and
   symplectic forms}, Math. Research Letters, 1 (1994) 809-822.


\noindent [T2] C. Taubes, The Seiberg-Witten invariants and
the Gromov invariants, Math. Research Letters 2 (1995), 221-238.


\noindent {[W]} C. T. C. Wall, { Diffeomorphisms of 4-manifolds},
      J. London  Math. Soc. 39 (1964) 131-140;


\bigskip
\noindent  Academy of Mathematics and Systems Science, Academia Sinica,
Beijing, 100080, P.R. China

\noindent libh@iss06.iss.ac.cn

\noindent Department of Mathematics, Princeton University, Princeton, NJ
08544, USA

\noindent tli@math.princeton.edu

\end